\documentclass[titlepage,11pt]{article}
\usepackage{latexsym}
\usepackage{amsfonts}
\usepackage{amsmath}
\usepackage[T1]{fontenc}
\newcommand\blackslug{\hbox{\hskip 1pt \vrule width 4pt height 8pt depth 1.5pt
        \hskip 1pt}}
\newcommand\bbox{\hfill \quad \blackslug \bigbreak}

\def\ll{,\ldots,}
\def\cupcup{\cup\cdots\cup}

\title{Polynomial bounds for chromatic number.\\ III. Excluding a double star}
\author{Alex Scott\thanks{Research supported by EPSRC grant EP/V007327/1.}\\
Mathematical Institute, University of Oxford, Oxford OX2 6GG, UK
\\
\\
Paul Seymour\thanks{Supported by AFOSR grant
A9550-19-1-0187, and by NSF grant  DMS-1800053.}\\
Princeton University, Princeton, NJ 08544
\\
\\
Sophie Spirkl\thanks{We acknowledge the support of the Natural Sciences and Engineering Research
Council of Canada (NSERC) [funding reference number RGPIN-2020-03912].
Cette recherche a \'et\'e financ\'ee par le Conseil de recherches en sciences
naturelles et en g\'enie du Canada (CRSNG) [num\'ero de r\'ef\'erence
RGPIN-2020-03912].}\\
University of Waterloo, Waterloo, Ontario N2L3G1, Canada}

\date{}

\newtheorem{thm}{}[section]

\newcommand{\Proof}{\noindent{\bf Proof.}\ \ }

\begin{document}
\maketitle
\begin{abstract}
A ``double star'' is a tree with two internal vertices. It is known that the Gy\'arf\'as-Sumner conjecture holds for double stars,
that is, for every double star $H$, there
is a function $f_H$ such that if $G$ does not contain $H$ as an induced subgraph then $\chi(G)\le f_H(\omega(G))$ (where $\chi, \omega$ are the chromatic number
and the clique number of $G$). Here we prove that $f_H$ can be chosen to be a polynomial.
\end{abstract}

\section{Introduction}
A class of graphs is {\em hereditary} if it is closed under isomorphism and under taking induced subgraphs. 
A hereditary class of graphs $\mathcal C$ is {\em $\chi$-bounded} if there is a function $f$ such that $\chi(G) \leq f(\omega(G))$ for every graph
$G\in\mathcal C$, where $\chi(G)$ and $\omega(G)$ denote the chromatic number and the clique number of $G$. There is a large literature addressing the question of which graph classes are $\chi$-bounded, and many open questions (see~\cite{survey} for a survey).

Hereditary classes defined by excluding some fixed graph $H$ are of particular interest. 
If $G,H$ are graphs, we say $G$ is {\em $H$-free} if no induced subgraph of $G$ is isomorphic
to $H$.
It is easily 
seen that if the class of $H$-free graphs is $\chi$-bounded then $H$ must be a forest, as Erd\H{o}s~\cite{erdos} showed that there are graphs with arbitrarily large girth and 
chromatic number.  The famous Gy\'arf\'as-Sumner conjecture~\cite{gyarfas, sumner} asserts the converse:
\begin{thm}\label{GSconj}
{\bf Conjecture: } For every forest $H$, there is a function $f$
such that $\chi(G)\le f(\omega(G))$ for every $H$-free graph $G$. 
\end{thm}
The Gy\'arf\'as-Sumner conjecture remains open in general, though it has been proved for some 
very restricted families of trees
(see, for example,~\cite{distantstars, gyarfasprob, gst, kierstead, kierstead2, scott, newbrooms, survey}).
In particular, it was proved by Kierstead and Penrice~\cite{kierstead} for trees of radius two.

Louis Esperet~\cite{esperet} made the striking conjecture that, for every $\chi$-bounded class $\mathcal C$, the function $f$ can be chosen to 
be a polynomial  (see the survey by Schiermeyer and Randerath~\cite{Schiermeyer}
for results on polynomial $\chi$-boundedness). In particular, this would imply a strengthening of the Gy\'arf\'as-Sumner conjecture, that
the function $f$ in \ref{GSconj} can always be chosen to be a polynomial. 
This is a bold conjecture, as frequently, when classes are known to be $\chi$-bounded, the best known function $f$ grows quite rapidly, often 
because the proofs use multiple applications of Ramsey-type results. Nevertheless, Esperet's strengthening has been verified for some cases of the Gy\'arf\'as-Sumner 
conjecture, for instance
when $H$ is a star, or a four-vertex path, or a matching (see~\cite{Schiermeyer}); and recently it has been shown when $H$ is obtained from a star 
by subdividing one edge once~\cite{liu},
 and when $H$ is a forest of stars~\cite{polystar}.

A {\em double star} is a tree in which at most two vertices have degree more than one.  Double stars have radius at most two, and so the result 
of Kierstead and Penrice~\cite{kierstead} shows that the class of $H$-free graphs is $\chi$-bounded whenever $H$ is a double star.  
In this paper, we prove a polynomial bound. 
Our main result is:
\begin{thm}\label{mainthm}
For every double star $H$, there is a polynomial $f$ such that $\chi(G)\le f(\omega(G))$ for every $H$-free graph $G$.
\end{thm}

This extends a theorem of Liu, Schroeder, Wang and Yu~\cite{liu}, who proved the same for double stars $H$ that
have at most one vertex of degree more than two. 

Our result is partially motivated by the Erd\H{o}s-Hajnal conjecture. In view of the recent result~\cite{fivehole}, the five-vertex path
$P_5$ is the smallest open case of this conjecture.  It is known that $P_5$ 
satisfies the Gy\'arf\'as-Sumner conjecture (in fact the 
Gy\'arf\'as-Sumner conjecture holds whever $H$ is a path), and if $P_5$ satisfies Esperet's strengthening then $P_5$ also satisfies the 
Erd\H{o}s-Hajnal conjecture. Thus $P_5$ appears likely to be a sticking point. We have not settled that; but this paper proves 
Esperet's strengthening for all trees that do not contain $P_5$.

We use standard notation throughout.
When $X\subseteq V(G)$, $G[X]$ denotes the subgraph induced on $X$. We write $\chi(X)$ for $\chi(G[X])$, and 
$\omega$ for $\omega(G)$, when there is no ambiguity.

\section{A degeneracy variant of defective colouring}
A graph $G$ is {\em $d$-degenerate}, or has {\em degeneracy at most $d$}, if every non-null subgraph 
$H$ has a vertex with degree (in $H$) at most $d$. Every $d$-degenerate graph has chromatic number 
at most $d+1$. 

Let us say a {\em $(k,d)$-colouring} of a graph is a partition $(A_1\ll A_k)$ of the vertex set $V(G)$, such that for $1\le i\le k$,
the subgraph induced on $A_i$ has degeneracy at most $d$; and we say that $G$ or $V(G)$ is {\em $(k,d)$-colourable} if there is such a partition. Thus, if a graph is 
$(k,d)$-colourable, its chromatic number is at 
most $k(d+1)$. We call this ``degenerate colouring''; it is a relative of ``defective colouring'', where we ask that the subgraph induced on each $A_i$
has maximum degree at most $d$, but it is not exactly the same (see~\cite{wood} for a survey of defective colouring). Let us explain
why we need to use degenerate colourings.

A standard way to bound the chromatic number of a graph $G$ is to partition $V(G)$ into some number of parts $V_1,V_2,\dots$, 
and bound the chromatic numbers of the parts
separately, and add to get a bound on $\chi(G)$. But we will be trying to prove that $\chi(G)$ is at most $\omega(G)^d$, for
some appropriately large constant $d$.
So for this ``addition'' method to work when $\omega$ is large,
if the best bound we know for one of the parts is something like $(\omega(G)-1)^d$, we would need much better bounds for 
all the other parts.

Fix a double star $H$, and choose a large constant $d$; and suppose that we
try to prove by induction on $\omega(G)$
that every $H$-free graph $G$ has chromatic number at most $\omega(G)^d$. The proof (by our method) does not work.
There comes a stage where $V(G)$ is partitioned into an unbounded number of parts $V_1,V_2,\ldots$. We will know, from
the induction on $\omega(G)$, that
each part has chromatic number at most something like $(\omega-1)^d$ (where $\omega=\omega(G)$), but we will not know a better bound for {\em any} of the parts. 
The ``addition'' method given above will therefore fail miserably. But we will know something about the edges between parts, 
which we
might hope will save us (though in fact it will not).  We will know that
for each part $V_i$, each of its vertices has only a small number of neighbours in the 
union of the later parts $V_{i+1}\cupcup V_n$;
say at most
$\omega^r$ neighbours, where $r$ is much less than $d$. Of course if there were no edges between the parts, all would be fine, and 
one might hope that similarly, because of the sparseness of the edges between parts, the effect of these edges could be fitted into the 
difference between
$\omega^d$ and $(\omega-1)^d$. But we can't do this (at least with no further information); the effect is multiplicative rather than 
additive. Even if for each $i$, each vertex of $V_i$ has only at most one neighbour in $V_{i+1}\cupcup V_n$, the chromatic number of the 
union of the parts might be $3/2$ times the maximum chromatic number of the individual parts, which is much too big.

Thus the inductive proof that every $H$-free graph $G$ has chromatic number at most $\omega(G)^d$ fails; and
for that reason we will instead prove by induction a stronger statement, about degenerate colourings.
We will prove by induction on $\omega$ that if $G$ is $H$-free, then 
$G$ is $(\omega^d, \omega^{r+1})$-colourable.
Then, when the situation above arises, we will know that each part
admits an $((\omega-1)^d, (\omega-1)^{r+1})$-colouring 
and hence an $(\omega^d, (\omega-1)\omega^{r})$-colouring. The union
of these colourings becomes an $(\omega^d, \omega^{r+1})$-colouring, which is what we want, and now it all works. 
Induction on $\omega(G)$ will be used to prove the statement about degenerate colouring, and then we 
deduce the statement about normal colouring at the end.

Let us state formally the lemma we just mentioned:
\begin{thm}\label{chaining}
Let $k,d,d'\ge 0$ be integers. Let $V(G)$ be partitioned into $V_1,V_2\ll V_n$, such that
\begin{itemize}
\item for $1\le i\le n$, $G[V_i]$ admits a $(k,d)$-colouring; and
\item for $1\le i<n$, every vertex of $V_i$ has at most $d'$ neighbours in $V_{i+1}\cupcup V_n$.
\end{itemize}
Then $G$ admits a $(k,d+d')$-colouring.
\end{thm}
The proof is clear.

\section{Templates}

The paper by Kierstead and Penrice~\cite{kierstead} uses the method of ``templates'', an idea that was introduced 
in~\cite{gst}  and has been applied in several papers to prove special cases of \ref{GSconj}. We will use the same idea,
but substantially modified to keep the numbers polynomial. Fix an integer $s$. (Eventually $s+1$ will be the maximum degree in the 
double star we are excluding.) 
We say an {\em $s$-template} in a graph $G$ is
a sequence $\mathcal{L}$ of pairwise disjoint subsets $(L_0,L_1\ll L_k)$, with $k\ge 0$, 
such that:
\begin{itemize}
\item  $L_0$ is a clique of $G$ (possibly empty), and every vertex in $L_0$ is adjacent to every vertex in $L_1\cupcup L_k$; 
\item $\omega^{s+5}\le |L_i|\le 14\omega^{s+6}$ for $1\le i\le k$; and
\item for all distinct $i,j\in \{1\ll k\}$, each vertex in $L_i$ has at most $\omega^{s+3}$ non-neighbours in $L_j$.
\end{itemize}
We say that $k$ is the {\em length} of the $s$-template, and define its 
{\em value} to be
$$|L_1\cupcup L_k|+7\omega^{s+5}|L_0|+k\omega^{s+5}.$$
Let us define $V(\mathcal{L})=L_0\cup L_1\cupcup L_k$, and $N(\mathcal{L})$ to be the 
set of vertices in $V(G)\setminus V(\mathcal{L})$ that have a neighbour in $L_1\cupcup L_k$. (Note that we do not consider 
neighbours in $L_0$.)

The idea of the proof is as follows. Let $G$ be an $H$-free graph (where $H$ is a double star).
Suppose that $G$ has chromatic number at least some huge (but some fixed constant) power of 
$\omega$. It follows from a theorem of an earlier paper~\cite{Kst} that $G$ has a subgraph
which is a complete bipartite graph, in which both parts have cardinality $14\omega^{s+6}$.
Consequently $G$ contains an $s$-template of length two and value at least $28\omega^{s+6}$. 
Thus we may choose
an $s$-template $\mathcal{L}$ with maximum value, and its value will also be at least $28\omega^{s+6}$. We have tuned these values
so that such $s$-templates have many useful properties.
The bulk of the proof is to bound the chromatic number of the set $N(\mathcal{L})$ (more exactly, to show it admits a certain 
degenerate colouring). Having done that, let $Z$ consist of $L_0$ together with all vertices in $N(\mathcal{L})$ that have at least a few 
neighbours that are not in $N(\mathcal{L})$ 
(a ``few'' means a constant power of $\omega$). It is easy to show that $|Z|$ is at most another constant power of $\omega$, and so 
every vertex of $(V(\mathcal{L})\cup N(\mathcal{L}))\setminus Z$ has only a few neighbours in the complement of this set.
It remains to bound the 
chromatic number of the complementary set, that is 
of $(V(G)\cup Z)\setminus (V(\mathcal{L})\cup N(\mathcal{L}))$; and to do this, we choose an $s$-template in this graph
with value as large as possible, and choose a third
with no neighbours in the first or second, and so on. This is the situation we discussed in the previous section, which motivated us
to use degenerate colouring, and this will allow us to produce a degenerate colouring of the whole of $G$.

We begin in this section by proving some properties of {\em optimal} $s$-templates, that is, $s$-templates chosen with maximum value.

\begin{thm}\label{clique}
Let $(L_0,L_1\ll L_k)$ be an $s$-template in $G$. There is a clique of $G$ with one vertex in each
of $L_1\ll L_k$, and consequently $k+|L_0|\le \omega$.
\end{thm}
\Proof Choose $v_1\in L_1$; and inductively for $2\le i\le k$, having chosen $v_1\ll v_{i-1}$, choose $v_i\in L_i$ as follows. 
There are at most $\omega^{s+3}$ vertices in $L_i$ nonadjacent to $v_h$, for $1\le h<i$; and since $\{v_1\ll v_{i-1}\}$ is a clique 
and therefore $i-1\le \omega$, it follows that $\omega^{s+3}(i-1)\le \omega^{s+4}< \omega^{s+5}\le |L_i|$. Consequently there exists
$v_i\in L_i$ adjacent to all of $v_1\ll v_{i-1}$. This completes the inductive definition of $v_1\ll v_k$. Hence 
$L_0\cup \{v_1\ll v_k\}$ is a clique, and so $k+|L_0|\le \omega$. This proves \ref{clique}.~\bbox

\begin{thm}\label{nottoosmall}
Let $(L_0,L_1\ll L_k)$ be an optimal $s$-template in $G$. If $\mathcal{L}$ has value at 
least $28\omega^{s+6}$, then $k\ge 2$, and 
$|L_i|\ge 5\omega^{s+5}$ for $1\le i\le k$.
\end{thm}
\Proof 
If $k=0$, the $s$-template has value $7\omega^{s+5}|L_0|\le 7\omega^{s+6}<28\omega^{s+6}$, a contradiction. If $k=1$, then since
$|L_1|\le 14\omega^{s+6}$, and $|L_0|+k\le \omega$ by \ref{clique}, the $s$-template has value at most 
$$14\omega^{s+6}+7\omega^{s+5}|L_0|+\omega^{s+5}\le 14\omega^{s+6}+7\omega^{s+5}(\omega-1)+\omega^{s+5}\le 21\omega^{s+6}<28\omega^{s+6},$$
a contradiction. So $k\ge 2$.

By reordering $L_1\ll L_k$ we may assume that $|L_1|\ll |L_h|\ge  5\omega^{s+5}$ and $|L_{h+1}|\ll |L_k|< 5\omega^{s+5}$. 
We will show that $h=k$.
For $h+1\le i\le k$, choose $v_i\in L_i$ such that $\{v_{h+1}\ll v_k\}$ is a clique $X$ (this is possible 
by \ref{clique}). For $1\le i\le h$, let $L_i'$ be the set of all vertices in $L_i$ that are adjacent to every vertex of $X$.
Since each vertex of $X$ has at most $\omega^{s+3}$ non-neighbours in $L_i$, it follows that 
$$|L_i'|\ge |L_i|-\omega^{s+4}\ge 5\omega^{s+5}-\omega^{s+5}\ge \omega^{s+5}$$
for $1\le i\le h$. 
Consequently 
$$(L_0\cup \{x_{h+1}\ll x_k\}, L_1'\ll L_h')$$ 
is an
$s$-template in $G$. Its value is that of $(L_0, L_1\ll L_k)$ plus 
$$7\omega^{s+5}(k-h)-\omega^{s+5}(k-h)-\left(|L_1\cupcup L_k|-|L_1'\cupcup L_h'|\right),$$
and so this is at most zero, since $(L_0\ll L_k)$ is optimal.
But 
$$|L_1\cupcup L_k|-|L_1'\cupcup L_h'|\le \sum_{1\le i\le h}\left(|L_i|-|L_i'|\right) +\sum_{h+1\le i\le k}|L_i|\le h\omega^{s+4}+
(k-h)(5\omega^{s+5}-1),$$
and consequently
$$7\omega^{s+5}(k-h)-\omega^{s+5}(k-h) \le h\omega^{s+4}+(k-h)(5\omega^{s+5}-1),$$ 
that is, 
$$(\omega^{s+5}+1)(k-h) \le h\omega^{s+4}.$$
But $\omega^{s+5} \ge  h\omega^{s+4},$
and so 
$(\omega^{s+5}+1)(k-h)\le omega^{s+5}$, which implies that $h=k$.
This proves \ref{nottoosmall}.~\bbox

\begin{thm}\label{manynonnbrs}
Let $(L_0,L_1\ll L_k)$ be an optimal $s$-template in $G$. If its value is at
least $28\omega^{s+6}$, 
then for $1\le i\le k$, every vertex in $L_i$ has at least $4\omega^{s+5}$ non-neighbours in $L_i$.
\end{thm}
\Proof
Let $i = 1$ say, and let $v\in L_1$. Let $L_1'$ be the set of neighbours of $v$ in $L_1$, and $M=L_1\setminus (L_1'\cup \{v\})$. 
We will show that $|M|\ge 4\omega^{s+5}$. If $|L_1'|<\omega^{s+5}$, then
$v$ has at least $|L_1|-\omega^{s+5}$ non-neighbours in $L_1$; and since $|L_1|\ge 5\omega^{s+5}$ by \ref{nottoosmall}, it follows that
$|M|\ge 4\omega^{s+5}$ as required. Thus we may assume that $|L_1'|\ge \omega^{s+5}$. 
For $2\le i\le k$ let $L_i'$ be the set of vertices in $L_i$ adjacent to $v$; thus by \ref{nottoosmall}, 
$$|L_i'|\ge |L_i|-\omega^{s+3}\ge 5\omega^{s+5}-\omega^{s+3}\ge \omega^{s+5}$$
for $2\le i\le k$. Consequently 
$$(L_0\cup \{v\}, L_1', L_2'\ll L_k')$$
is an $s$-template. From the optimality of $(L_0\ll L_k)$, it follows that
$7\omega^{s+5}-(|M|+1)-(k-1)\omega^{s+3}\le 0$, and so 
$$|M|\ge 7\omega^{s+5}-1-(k-1)\omega^{s+3}\ge 4\omega^{s+5}.$$ 
This proves \ref{manynonnbrs}.~\bbox

\begin{thm}\label{linked}
Let $(L_0,L_1\ll L_k)$ be an optimal $s$-template in $G$.
Suppose that $\omega\ge 4$, and 
$A,B$ are disjoint subsets of $L_1$, with $|L_1\setminus (A\cup B)|\le \omega^{s+3}$, such that every vertex in $A$ has fewer than
$\omega^s$ non-neighbours in $B$. Then either $|B|<14\omega^{s+1}$ or $A=\emptyset$.
\end{thm}
\Proof
Suppose that $|B|\ge 14\omega^{s+1}$.
\\
\\
(1) {\em $|B|\ge 2\omega^{s+5}$.}
\\
\\
Suppose that
$|B|<2\omega^{s+5}$. Each vertex in $B$ has at least $4\omega^{s+5}$ non-neighbours in $L_1$, by \ref{manynonnbrs}, and only
at most $2\omega^{s+5}+\omega^{s+3}$ of them do not belong to $A$; and since $2\omega^{s+5}-\omega^{s+3}\ge \omega^{s+5}$, there are at least $\omega^{s+5}|B|\ge 14\omega^{2s+6}$
nonedges between $B$ and $A$. Since $|A|\le |L_1|\le 14\omega^{s+6}$, some vertex in $A$
has at least $\omega^s$ non-neighbours in $B$, a contradiction. This proves (1).
\\
\\
(2) {\em $|A|<\omega^{s+5}$.}
\\
\\
Suppose that $|A|\ge \omega^{s+5}$. Let $B'$ be the set of all vertices in $B$ with at most $\omega^{s+3}$ non-neighbours in $A$. Since there
are only at most $\omega^s|A|$ nonedges between $A,B$, there are at most $|A|/\omega^{3}$ vertices in $B$ that have more than $\omega^{s+3}$
non-neighbours in $A$; and so $|B'|\ge |B|-|A|/\omega^{3}\ge \omega^{s+5}$, by (1) and since $|A|/\omega^{3}\le 14\omega^{s+3}< \omega^{s+5}$ (the last because $\omega\ge 4$).
Hence
$$(L_0, A,B', L_2\ll L_k)$$
is an $s$-template, and the optimality of $(L_0\ll L_k)$ implies that
$$|A|+|B'|+\omega^{s+5}\le |L_1|\le |A|+|B'|+\omega^{s+3}+|A|/\omega^{3},$$
and so
$\omega^{s+5}\le \omega^{s+3}+14\omega^{s+3},$ a contradiction. This proves (2).

\bigskip
Suppose that $A\ne \emptyset$, and choose $v\in A$. Since $v$ has at least $4\omega^{s+5}$ non-neighbours in $L_1$ by \ref{manynonnbrs},
and at most $\omega^{s+5}$ of them belong to $A$ by (2), and at most $\omega^{s+3}$ are not in $A\cup B$, it follows that $v$ has at least
$3\omega^{s+5}-\omega^{s+3}\ge \omega^s$
non-neighbours in $B$, a contradiction. Thus $A=\emptyset$. This proves \ref{linked}.~\bbox

\begin{thm}\label{dense}
Let $\mathcal{L}=(L_0,L_1\ll L_k)$ be an optimal $s$-template in $G$.
There are fewer than $14\omega^{s+6}$ vertices $v\in N(\mathcal{L})$ such that for $1\le i\le k$, 
$v$ has at most $\omega^{s+2}/4$
non-neighbours in $L_i$.
\end{thm}
\Proof
Suppose that there is a set $M\subseteq N(\mathcal{L})$ with $|M|=14\omega^{s+6}$, such that 
for each $v\in M$ and for $1\le i\le k$, $v$ has at most $\omega^{s+2}/4$
non-neighbours in $L_i$. For $1\le i\le k$, there are at most $|M|\omega^{s+2}/4=7\omega^{2s+8}/2$ nonedges between $M$ and $L_i$; and so 
there are
at most $7\omega^{s+5}/2$ vertices in $L_i$ with at least $\omega^{s+3}$ non-neighbours in $M$. Let $L_i'$ be the set of vertices in $L_i$ that have
fewer than $\omega^{s+3}$ non-neighbours in $M$. Hence $|L_i'|\ge |L_i|-(7\omega^{s+5}/2)\ge \omega^{s+5}$, since $|L_i|\ge 5\omega^{s+5}$ by \ref{nottoosmall}.
It follows that 
$$(\emptyset, L_1'\ll L_k',M)$$
is an $s$-template. Its value is 
$$|L_1'\cupcup L_k'|+|M|+(k+1)\omega^{s+5}\ge |L_1\cupcup L_k|-k(7\omega^{s+5}/2)+14\omega^{s+6}+(k+1)\omega^{s+5}$$
and the optimality of $(L_0\ll L_k)$ implies that
$$|L_1\cupcup L_k|-k(7\omega^{s+5}/2)+14\omega^{s+6}+(k+1)\omega^{s+5}\le |L_1\cupcup L_k|+7\omega^{s+5}|L_0|+k\omega^{s+5},$$
that is, 
$$2\omega+1/7\le |L_0|+k/2,$$
contrary to \ref{clique}. This proves \ref{dense}.~\bbox

\section{Using the double star}
We are concerned with graphs that do not contain some fixed double star, but so far we have not used that fact. For $s\ge 1$, 
let $H_s$
be the double star with $2s+2$ vertices, with two internal vertices both of degree $s+1$. Every double star is an induced subgraph of $H_s$
for some $s$, so it suffices to prove the result for $H_s$-free graphs. 

We will need the following result of~\cite{Kst}:
\begin{thm}\label{Kst}
Let $H$ be a forest. Then there exists $c>0$ such that for every $H$-free graph $G$ and every integer $t\ge 0$, either
$G$ contains the complete bipartite graph $K_{t,t}$ as a subgraph, or $G$ has degeneracy less than $t^c$, and hence has
chromatic number at most $t^c$.
\end{thm} 

We will also need the following version of Ramsey's theorem (well-known, but proved for instance in~\cite{polystar}):
\begin{thm}\label{Ramsey}
If $s\ge 0$ is an integer, then every graph $G$ with no stable set of cardinality $s$ has at most
$$\omega^{s-1}+\omega^{s-2}+\cdots+\omega$$
vertices, and hence fewer than $\omega^s$ vertices if $\omega>1$.
\end{thm}

Let $\mathcal{L}=(L_0\ll L_k)$ be an optimal $s$-template in $G$. With respect to this template, we say a vertex 
$v\in N(\mathcal{L})$ is 
\begin{itemize}
\item {\em pendant}
if there exist distinct $i,j\in \{1\ll k\}$, and a vertex $u\in L_j$, and a stable set $S$ of
$s+1$ vertices in $L_i$, all adjacent to $u$, such that $v$ is not adjacent to $u$, and $v$ has exactly one neighbour in $S$;
\item 
{\em dense}
if there exists $j\in \{1\ll k\}$, and $u\in L_j$, such that for all
$i\in \{1\ll k\}\setminus \{j\}$, there are fewer than $\omega^{s+2}/14$ vertices in $L_i$ that are adjacent to $u$ and not to $v$;
\item {\em pure} if there are least two values of $i\in \{1\ll k\}$
such that $v$ has no neighbour in $L_i$, and for $1\le i\le k$, 
either $v$ has no neighbour in $L_i$ or $v$ has at most $\omega^{s+2}/7$ non-neighbours in $L_i$.
\end{itemize}

\begin{thm}\label{types}
Let $s\ge 1$ be an integer, and let $G$ be $H_s$-free, with $\omega\ge 200$. Let $\mathcal{L}=(L_0\ll L_k)$ be an optimal $s$-template in $G$.
Then every vertex in $N(\mathcal{L})$ is either pendant or dense or pure with respect to $\mathcal{L}$.
\end{thm}
\Proof
Let $v\in N(\mathcal{L})$, and suppose that $v$ is neither dense nor pendant with respect to $\mathcal{L}$. We will prove that $v$ 
is pure. 
\\
\\
(1) {\em For all distinct $i,j\in \{1\ll k\}$, if $u\in L_j$ is nonadjacent to $v$, then either $u,v$ have no common neighbour 
in $L_i$, or fewer than $\omega^{s+2}/14$ neighbours of $u$ in $L_i$ are nonadjacent to $v$.}
\\
\\
Let $A$ be the set of all vertices
in $L_i$ adjacent to both $u,v$, and let $B$ be the set of vertices in $L_i$ adjacent to $u$ and not to $v$. 
Suppose that $A\ne \emptyset$, and $|B|\ge \omega^{s+2}/14$.
Since  $|L_i\setminus (A\cup B)|\le \omega^{s+3}$ (because $u$ has at most $\omega^{s+3}$ non-neighbours in $L_i$),
and $|B|\ge \omega^{s+2}/14\ge 14\omega^{s+1}$ (because $\omega\ge 200$), \ref{linked} implies that some vertex $w\in A$ has at least 
$\omega^s$ non-neighbours in $B$. By \ref{Ramsey}, this set of non-neighbours includes a stable set of size $s$, contradicting that $v$
is not pendant. Thus either $A=\emptyset$ or $|B|<\omega^{s+2}/14$. This proves (1).
\\
\\
(2) {\em For all $j\in \{1\ll k\}$, if $v$ has a non-neighbour in $L_j$, there exists $i\in \{1\ll k\}$ with $i\ne j$
such that $v$ has at most $\omega^{s+3}$ neighbours in $L_i$.}
\\
\\
Choose $u\in L_j$ nonadjacent to $v$. Since $v$ is not dense, there exists $i\in \{1\ll k\}$ with $i\ne j$
such that there are at least $\omega^{s+2}/14$ vertices in $L_i$ adjacent to $u$ and not to $v$. By (1), $u,v$ have
no common neighbour in $L_i$, and hence $v$ has at most $\omega^{s+3}$ neighbours in $L_i$.
This proves (2).
\\
\\
(3) {\em For $1\le j\le k$, if $v$ has at most $\omega^{s+3}$ neighbours in $L_j$, then $v$ has no neighbours in $L_j$.}
\\
\\
By (2), there exists $i\in \{1\ll k\}$ different from $j$, 
such that $v$ has at most $\omega^{s+3}$ neighbours in $L_i$. Suppose that $v$ has a neighbour $w\in L_j$. Since $w$ has at most $\omega^{s+3}$
non-neighbours in $L_i$, and $v$ has at most $\omega^{s+3}$ neighbours in $L_i$, and $|L_i|\ge \omega^{s+5}>2\omega^{s+3}$, there exists $x\in L_i$
adjacent to $w$ and not to $v$. Now $w$ has at least $4\omega^{s+5}$ non-neighbours in $L_j$, by \ref{manynonnbrs}; at most $\omega^{s+3}$
of them are nonadjacent to $x$, and at most $\omega^{s+3}$ of them are adjacent to $v$, and so at least $4\omega^{s+5}-2\omega^{s+3}\ge \omega^s$
of them are nonadjacent to $v$ and adjacent to $x$. By \ref{Ramsey}, this set includes a stable set of size $s$, 
and so $v$ is pendant, a contradiction.  Thus $v$ has no neighbour in $L_j$. This proves (3).
\\
\\
(4) {\em For $1\le i\le k$, either $v$ has at most $\omega^{s+2}/7$ non-neighbours in $L_i$, or $v$ has no neighbours in $L_i$.}
\\
\\
Let $A_i,B_i$ be the sets of neighbours
and non-neighbours respectively of $v$ in $L_i$, and suppose that $A_i\ne \emptyset$, and $|B_i|>\omega^{s+2}/7$.
By (2) and (3), there
exists $j\in \{1\ll k\}$ with $j\ne i$ such that $v$ has no neighbours in $L_j$. 
By (1), for each $u\in L_j$, either $u$ has no neighbours in $A_i$, or $u$ has at most $\omega^{s+2}/14$ neighbours in $B_i$.
Let $X$ be the set of vertices in $L_j$ with no neighbour in $A_i$, and $Y=L_j\setminus X$. Since $A_i\ne \emptyset$, and a
vertex in $A_i$ has at most $\omega^{s+3}$ non-neighbours in $L_j$, it follows that $|X|\le \omega^{s+3}$, and so $|Y|\ge |L_j|-\omega^{s+3}$. Every vertex in $Y$
is adjacent to at most half the vertices in $B_i$ (since $|B_i|\ge \omega^{s+2}/7$), and so some vertex $b\in B_i$
is adjacent to at most half the vertices in $Y$. Since $b$ has at most $\omega^{s+3}$ non-neighbours in $L_j$, it follows that
$|Y|/2\le \omega^{s+3}$; but $|X|\le \omega^{s+3}$, and so $|L_j|\le 3\omega^{s+3}$, a contradiction. This proves (4).

\bigskip

Since $v$ is not dense, it has a non-neighbour in one of $L_1\ll L_k$; and so by (2) and (3), it has no neighbours in some
$L_j$. By (2), $v$ has at most $\omega^{s+3}$ neighbours in $L_i$ for some $i\ne j$; and so has no neighbours in $L_i$
by (3). From (4) it follows that $v$ is pure. This proves \ref{types}.~\bbox

\begin{thm}\label{pendant}
Let $s\ge 1$ be an integer, and let $G$ be $H_s$-free. Let $(L_0\ll L_k)$ be an optimal $s$-template in $G$.
There are at most $14^{s+2}\omega^{s^2+9s+14}$ pendant vertices.
\end{thm}
\Proof
Let $i,j\in \{1\ll k\}$ be distinct, let $u\in L_j$, let $S\subseteq L_i$ be a stable set of $s+1$ neighbours of $u$, and let
$w\in S$. Let $X(i,j,u,S,w)$ be the set of all $v\in V(G)\setminus (L_0\cupcup L_k)$ such that $v$ is adjacent to $w$
and is nonadjacent to all other vertices in $S\cup \{u\}$.
If $X(i,j,u,S,w)$ includes a stable set of size $s$, say $T$, then the subgraph
induces on $S\cup T\cup \{u\}$ is isomorphic to $H_s$, a contradiction. Thus, \ref{Ramsey} implies that $|X(i,j,u,S,w)|\le \omega^s$.
Since there are only $k^2(14\omega^{s+6})^{s+2}$ choices for $i,j,u,S,w$, and every pendant vertex belongs to $X(i,j,u,S,w)$
for some choice of $i,j,u,S,w$, it follows that there are at most $k^2(14\omega^{s+6})^{s+2}\omega^s$ pendant
vertices. Since $k\le \omega$ by \ref{clique}, this proves \ref{pendant}.~\bbox

\bigskip

\begin{thm}\label{semidense}
Let $s\ge 1$ be an integer, and let $G$ be $H_s$-free. Let $\mathcal{L}=(L_0\ll L_k)$ be an optimal $s$-template in $G$.
Let $c$ satisfy \ref{Kst} when $H=H_s$.
The chromatic number of the set of all dense vertices is at most $(14\omega^{s+6})^{c+1}\omega$.
\end{thm}
\Proof Let $1\le j\le k$, let $u\in L_j$, and let $X(j,u)$ be the set of all $v\in N(\mathcal{L})$ such that
for all
$i\in \{1\ll k\}\setminus \{j\}$, there are fewer than $\omega^{s+2}/14$ vertices in $L_i$ that are adjacent to $u$ and not to $v$.

Suppose that $\chi(X(k,u))> (14\omega^{s+6})^c$ for some $u\in L_k$. Then by \ref{Kst}, $G[X(k,u)]$
contains a copy of $K_{14\omega^{s+6},14\omega^{s+6}}$ as a subgraph; let $M_1,M_2$ be disjoint subsets of $X(k,u)$, both of cardinality $14\omega^{s+6}$, such that
every vertex of $M_1$ is adjacent to every vertex of $M_2$. For $1\le i\le k-1$, let $L_i'$ be the set of vertices in $L_i\cap N(u)$
that have at most $\omega^{s+3}$ non-neighbours in $M_1$ and at most $\omega^{s+3}$ non-neighbours in $M_2$. There are at most 
$(\omega^{s+2}/14)14\omega^{s+6}$ nonedges
between $M_1$ and $L_i\cap N(u)$, and it follows that at most $\omega^{s+5}$ vertices in $L_i\cap N(u)$
have more than $\omega^{s+3}$ non-neighbours in $M_1$; and the same for $M_2$. Consequently
$$|L_i'|\ge |L_i\cap N(u)|-2\omega^{s+5}\ge 5\omega^{s+5} -\omega^{s+3}-2\omega^{s+5}\ge \omega^{s+5},$$
and so
$$(\emptyset,L_1'\ll L_{k-1}', M_1,M_2)$$
is an $s$-template. Moreover, $|L_i'|\ge |L_i|-\omega^{s+3}-2\omega^{s+5}$, since $u$ has at most $\omega^{s+3}$ non-neighbours in $L_i$.
From the optimality of $\mathcal{L}$, it follows that
$$28\omega^{s+6}\le |L_k|+(\omega^{s+3}+2\omega^{s+5})(k-1)+ 7\omega^{s+5}|L_0|\le 14\omega^{s+6}+7\omega^{s+6},$$ a contradiction. Consequently $\chi(X(k,u))\le (14\omega^{s+6})^c$. The union
of the sets $X(k,u)$ over all $u\in L_k$ thus has chromatic number at most $(14\omega^{s+6})^{c+1}$; and from the symmetry between $L_1\ll L_k$,
the set of all dense vertices has chromatic number at most $k(14\omega^{s+6})^{c+1}\le (14\omega^{s+6})^{c+1}\omega$. This proves \ref{semidense}.~\bbox

\section{Pure vertices}
In view of \ref{types}, \ref{pendant} and \ref{semidense}, in order to bound the chromatic number of $N(\mathcal{L})$ it remains to bound
the chromatic number of the set of pure vertices, and that is the topic of this section. 
But here we will need to use induction on $\omega$,
and so as we discussed earlier, we will in fact work with degenerate colouring. 

Throughout this section, let $s\ge 1$ 
be an integer, let $G$ be $H_s$-free, let $\mathcal{L}=(L_0\ll L_k)$ be an optimal $s$-template in $G$.
Let $M$
denote the set of all vertices that are pure with respect to this template. For each $v\in M$, let $I_v$ be the set of all 
$i\in \{1\ll k\}$ such that $v$ has a neighbour in $L_i$ (and hence has at most $\omega^{s+2}/14$ non-neighbours in $L_i$). 
For each $I\subseteq \{1\ll k\}$, let $M_I$ be the set of all $v\in M$ with $I_v=I$. 

We wish to find a degenerate colouring of the union of all the sets $M_I$. One problem is that the number of sets $I$
with $M_I$ nonempty may be superpolynomial; but we will show that there are only a linear number of sets $M_I$ of large 
cardinality, and we can colour all the small ones simultaneously.

\begin{thm}\label{incomp}
Let $I\subseteq \{1\ll k\}$ and $u\in M_I$ and $j\in \{1\ll k\}\setminus I$. 
For each $i\in I$, there are fewer than $\omega^{s}$ vertices $v$ adjacent to $u$ such that $v\in M_J$ for some 
$J\subseteq \{1\ll k\}\setminus \{i,j\}$.
Hence
there are fewer than $\omega^{s+1}$ vertices $v$ adjacent to $u$ such that $v\in M_J$ for some 
$J\subseteq \{1\ll k\}\setminus \{j\}$ with $I\not\subseteq J$.
\end{thm}
\Proof
Let $i\in I$, and let $\mathcal{A}_{i}$ be the set of all $J\subseteq \{1\ll k\}\setminus\{i,j\}$.
Let $V_{i}$ be the union of all the sets $M_J$ for $J\in \mathcal{A}_{i}$. Suppose that $u$ has $\omega^s$ neighbours in $V_{i}$.
By \ref{Ramsey}, there is a stable set $S\subseteq V_{i}$
of neighbours of $u$ with $|S|=s$. Choose $x\in L_i$ adjacent to $u$. Thus $x$ has no neighbour in $S$ since for each 
$J\in \mathcal{A}_{i}$, 
no vertex in $M_J$ has a neighbour in $L_i$. Since $x$ has at most $\omega^{s+3}$
non-neighbours in $L_j$, and $|L_j|\ge \omega^{s+5}\ge \omega^{s+3}+\omega^s$, it follows that $x$ has at least $\omega^s$ neighbours in $L_j$, and hence
by \ref{Ramsey} there is a stable set $T\subseteq L_j$ of neighbours of $x$ with $|T|=s$. Since $j\notin I\cup J$ for $J\in \mathcal{A}_{i}$, no vertex in 
$S\cup \{u\}$ has a neighbour in $T$; and so the subgraph induced on $S\cup T\cup \{u,x\}$ is isomorphic to $H_s$, a contradiction.

Hence $u$ has fewer than $\omega^s$ neighbours in $V_{i}$ for each choice of $i$. Since there are only at most $|I|\le \omega$
choices of $i$, 
this proves \ref{incomp}.~\bbox

For $m\ge 0$, we say $I\subseteq \{1\ll k\}$ 
is {\em $m$-small}
if $|M_I|\le m$, and {\em $m$-large} if it is not $m$-small.

\begin{thm}\label{smallchi}
For each $m\ge 0$, the union of the sets $M_I$ over all $m$-small $I\subseteq \{1\ll k\}$ has chromatic number at most 
$2\omega(m+\omega^{s+1})$.
\end{thm}
\Proof
For all $j\in \{1\ll k\}$, let $\mathcal{A}_{j}$ be the set of all $m$-small $I\subseteq \{1\ll k\}\setminus \{j\}$, and 
let $V_{j}$ be the union of the sets $M_I$ for all $I\in \mathcal{A}_{j}$. If $u,v\in V_j$ are adjacent, we will  direct the edge 
$uv$ as follows. Let $u\in M_I$ and $v\in M_J$ where $I,J\in \mathcal{A}_j$. If $|I|>|J|$ we direct the edge $uv$
from $u$ to $v$. If $|I|=|J|$ (and in particular, if $I=J$) we direct $uv$ arbitrarily. We claim every vertex $u\in V_j$
has outdegree less than $m+\omega^{s+1}$. Let $u\in M_I$. Certainly $u$ has at most $m$ out-neighbours in $M_I$ since $|M_I|\le m$.
If $v$ is an out-neighbour of $u$ and $v\in J\in \mathcal{A}_j$ where $J\ne I$, then $|J|\le |I|$, and
it follows that $I\not\subseteq J$, and $I\cup J\ne \{1\ll k\}$ since $j\notin I\cup J$; and so there are fewer than $\omega^{s+1}$
such vertices $v$ by \ref{incomp}. This proves that  every vertex $u\in V_j$
has outdegree less than $m+\omega^{s+1}$. Hence every subgraph of $G[V_j]$ with $n$ vertices has at most $n(m+\omega^{s+1}-1)$
edges, and so (if $n>0$) has a vertex of degree at most $2(m+\omega^{s+1}-1)$; and hence $G[V_j]$ has degeneracy at most 
$2(m+\omega^{s+1}-1)$, and so has chromatic number less than $2(m+\omega^{s+1})$. Since there are at most $\omega$ choices of $j$,
this proves \ref{smallchi}.~\bbox

It remains to handle the $m$-large sets.
For $t\ge 1$, let us say two disjoint subsets $A,B$ of $V(G)$ are {\em $s$-crowded} if there is no stable set $X$ with 
$|X\cap A|=|X\cap B|=s$.
\begin{thm}\label{nested}
Let $m\ge (s+1)\omega^s$, and let $I,J\subseteq \{1\ll k\}$ be $m$-large, with $I\ne J$. If $\omega^3>\omega+s/7$, then either:
\begin{itemize}
\item $J\subseteq I$ and every vertex in $M_I$ has fewer than $\omega^s$ neighbours in $M_J$; or
\item $I\subseteq J$ and every vertex in $M_J$ has fewer than $\omega^s$ neighbours in $M_I$; or
\item $I\cup J=\{1\ll k\}$ and $M_I, M_J$ are $s$-crowded.
\end{itemize}
Consequently there are at most $k-1$ $m$-large sets.
\end{thm}
\Proof
If $J\subseteq I$ then the first bullet holds by \ref{incomp}; so we may assume that $J\not\subseteq I$ and similarly $I\not\subseteq J$.
Choose $i\in I\setminus J$ and $j\in J\setminus I$. Suppose that $M_I, M_J$ are not $s$-crowded, and choose 
$S\subseteq M_I$ and $T\subseteq M_J$ with $|S|=|T|=s$ such that $S\cup T$ is stable.
Since each vertex in $S$ has at most $\omega^{s+2}/7$ non-neighbours in $i$, and $s\omega^{s+2}/7<\omega^{s+5}$, there exists $u\in L_i$
adjacent to every vertex in $S$. Since $u$ has at most $\omega^{s+3}$ non-neighbours in $L_j$, and each vertex in $T$
has at most $\omega^{s+2}/7$ non-neighbours in $L_j$, and $\omega^{s+3}+s\omega^{s+2}/7<\omega^{s+5}$, there exists $v\in L_j$ adjacent to $u$ and to each vertex in 
$T$. But then the subgraph induced on $S\cup T\cup \{u,v\}$ is isomorphic to $H_s$. This proves that $M_I, M_J$ are 
$s$-crowded.

Let $S\subseteq M_I$ be stable with $|S|=s$. (This is possible by \ref{Ramsey} since $m\ge \omega^s$.) Since $M_I, M_J$ are
$s$-crowded, \ref{Ramsey} implies that there are fewer than $\omega^s$ vertices in $M_J$ with no neighbour in $S$; 
and so some vertex in $S$
has at least 
$$(|M_J|-\omega^s)/s>(m-\omega^s)/s\ge \omega^s$$ 
neighbours in $M_J$. By \ref{incomp} it follows that $I\cup J=\{1\ll k\}$ and so the third bullet holds. This proves the
first assertion of \ref{nested}.

To show that there are at most $k-1$ $m$-large sets, for each $X\subseteq \{1\ll k\}$ with $|X|\le k-1$, 
let $n(X)$ be the number of $m$-large sets 
that include $X$. We prove by induction on $k-|X|$ that $n(X)\le k-|X|-1$. If $|X|=k-1$ then $n(X)=0$ and the claim holds,
so we assume that $|X|\le k-2$ and the result holds for all larger subsets of $\{1\ll k\}$. Since $k-|X|-1\ge 1$, we may assume that
at least two $m$-large sets include $X$, and so at least one properly includes $X$. Choose an $m$-large set $J$ minimal such 
that $X\subseteq J$ and $X\ne J$. If $I$ is an $m$-large set including $X$, then by \ref{incomp} and the minimality of $J$, 
either $I=X$, or $J\subseteq I$, or
$I\cup J=\{1\ll k\}$. Moreover, if $I\cup J=\{1\ll k\}$, then $I$ includes 
$K$, where $K=\{1\ll k\}\setminus (J\setminus X)$. Consequently 
$n(X)\le n(J)+n(K)+1$. Since $J,K$ are both proper supersets of $X$, and not equal to $\{1\ll k\}$,
the inductive hypothesis implies that $n(J)\le k-|J|-1$ and $n(K)\le k-|K|-1$, and so
$$n(X)\le (k-|J|-1)+(k-|K|-1)+1=k-|X|-1.$$
This completes the proof that $n(X)\le k-|X|-1$ for each $X\subseteq \{1\ll k\}$ with $|X|\le k-1$. Setting $X=\emptyset$, it follows
that there are only $k-1$ $m$-large sets. This proves \ref{nested}.~\bbox


\begin{thm}\label{largechi}
Suppose that $d\ge 1$ and $D\ge 0$ have the property that every $H_s$-free graph $G'$ with $\omega(G')<\omega(G)$ is $(\omega(G')^d,D)$-colourable, 
and let
$m= (s+1)\omega^s$. Assume that $\omega^3>\omega+s/7$
 and $\omega\ge 200$.
Then the union of the sets $M_I$ over all $m$-large $I$ is $(K,D)$-colourable where 
$$K=(\omega-1)^d+\omega^{s^2+4s+2}+\omega^{s+3}+1.$$
\end{thm}
\Proof Let $M^*$ be the union of the sets $M_I$ for all $m$-large $I$. We may assume that $M^*\ne \emptyset$.
\\
\\
(1) {\em There is a partition $\mathcal{A}, \mathcal{B}$ of the set of $m$-large sets, such that 
\begin{itemize}
\item $M_I, M_J$ are $s$-crowded for each $I\in \mathcal{A}$ and $J\in \mathcal{B}$; and
\item all the sets in $\mathcal{A}$ have an element in common, and so do all the sets in $\mathcal{B}$.
\end{itemize}
}
\noindent
There is an $m$-large set $J$ since $M^*\ne \emptyset$; choose $J$ minimal. Let
$\mathcal{A}$ be the set of all $m$-large sets that include $J$, and let $\mathcal{B}$ be the set of $m$-large sets that do not 
include $J$. Choose $j\in J$ and $j'\in \{1\ll k\}\setminus J$. ($J$ is nonempty since every pure vertex 
has a neighbour in $L_1\cupcup L_k$, and $J\ne \{1\ll k\}$ from the definition of pure.)
Every set in $\mathcal{A}$ includes $J$, and so contains $j$. We claim that every set in $\mathcal{B}$ includes 
$\{1\ll k\}\setminus J$ and so contains $j'$. To see this, let $I'\in \mathcal{B}$. Since $I'\in \mathcal{B}$, $I'$ does not 
include $J$. Also $J\not\supseteq I'$ from the minimality of $J$, and so from \ref{nested}, 
$I'\cup J=\{1\ll k\}$. This proves that every set in $\mathcal{B}$ includes
$\{1\ll k\}\setminus J$ and so contains $j'$. Finally, we must show that $M_I, M_{I'}$ are $s$-crowded for all 
$I\in \mathcal{A}$ and $I'\in \mathcal{B}$. Since $J\subseteq I$ and $J\not\subseteq I'$, it follows that $I\not\subseteq I'$;
and since $\{1\ll k\}\setminus J\subseteq I'$ and $\{1\ll k\}\setminus J\not\subseteq I$, it follows that $I'\not\subseteq I$.
By \ref{nested}, $M_I, M_{I'}$ are $s$-crowded. This proves (1).

\bigskip

Choose $\mathcal{A}, \mathcal{B}$ as in (1). Let $A$ be the union of all the sets $M_I$ for $I\in \mathcal{A}$, and 
define $B$ similarly. Thus $M^*=A\cup B$. 
\\
\\
(2) {\em Every clique included in $A$ has cardinality at most $\omega-1$, and the same for $B$.}
\\
\\
Suppose that there is a clique $X\subseteq A$ with
$|X|=\omega$. Let $j\in \{1\ll k\}$ belong to all the sets in $\mathcal{A}$. Since $(\omega^{s+2}/7)\omega<\omega^{s+5}\le |L_j|$, 
there exists $v\in L_j$ adjacent to every vertex of $X$, contradicting that $X$ is a clique of
$G$
of maximum cardinality. The same holds for $B$. This proves (2).

\bigskip
Let $n=\omega^{s+2}$. 
If $|A|\le n\omega$, then by (2) it follows that $G[M^*]$ admits an $(n\omega+(\omega-1)^d,D)$-colouring and the theorem holds. Thus we may assume that
 $|A|> n\omega$.
Let $X_1$ be the largest clique contained
in $A$, and inductively for $i\ge 2$ let $X_i$ be the largest
clique contained in $A\setminus (X_1\cupcup X_{i-1})$. Let $|X_n|=t$. Since $|A|> n\omega$ it follows that $t>0$.
Let $X=X_1\cupcup X_n$. Thus $|X|\le n\omega$, and $\omega(G[A\setminus X])\le t$.

Let $C$ be the set of all vertices $v\in B$ such that for some $I\in \mathcal{A}$, $v$ has at least $\omega^s$ non-neighbours in 
$M_I\cap X$. 
\\
\\
(3) {\em $|C|\le  n^s \omega^{2s+2}$.}
\\
\\
For each $v\in C$, there exists $I\in \mathcal{A}$ and a stable set $S\subseteq M_I\cap X$ with $|S|=s$ such that $v$
has no neighbour in $S$, by \ref{Ramsey}. For each such $I$ and $S$, and each $I'\in \mathcal{B}$, there are at most $\omega^s$ vertices in $M_{I'}$
that have no neighbours in $S$, by \ref{Ramsey} and since $M_I, M_{I'}$ are $s$-crowded (because $I\not\subseteq I'$ and 
$I'\not\subseteq I$ by (1), and by \ref{nested}). Consequently for each choice of $I, S$ there are at most $k\omega^s$
vertices in $C$ with no neighbour in $S$, since $|\mathcal{B}|\le k$ by \ref{nested}. For each choice of $I$ there
are only $|X|^s\le n^s\omega^s$ choices of $S$, since $|X|\le n\omega$; and there are only $k$ choices of $I$ by \ref{nested}.
Thus $|C|\le (k\omega^s)(n^s\omega^s)k\le n^s \omega^{2s+2}$. This proves (3).
\\
\\
(4) {\em Every clique in $B\setminus C$ has cardinality at most $\omega-t$.}
\\
\\
Let $Y\subseteq B\setminus C$ be a clique. Since $Y\cap C=\emptyset$, every vertex in $Y$
has at most $k\omega^s$ non-neighbours in $X$, since it has at most $\omega^s$ in each $X\cap M_I$ and there are only $k$
choices of $M_I$. Consequently at most $k\omega^s|Y|$ vertices in $X$ have a non-neighbour in $Y$. Since $n>k\omega^s|Y|$
(since $n=\omega^{s+2}$, and $k\le \omega$ and $|Y|<\omega$ by (2))
it follows that there exists $i\in \{1\ll n\}$ such that every vertex in $Y$ has no non-neighbour in $X_i$, and so $X_i\cup Y$
is a clique. But $|X|\ge t$ from the choice of $t$, and $|X\cup Y|\le \omega$, and so $|Y|\le \omega-t$. This proves (4).

\bigskip

Now $\chi(M^*)\le |X|+\chi(A\setminus X)+|C|+\chi(B\setminus C)$. But $|X|\le n\omega$; $A\setminus X$ is $(t^d, D)$-colourable, since 
$\omega(G[A\setminus X])\le t$ and $t<\omega$ (by (2)); $|C|\le n^s \omega^{2s+2}$ by (3); and $B\setminus C$ is $((\omega-t)^d,D)$-colourable
by (4) and since $\omega-t<\omega$ (because $t>0$). Thus $M^*$ is $(K_1,D)$-colourable where
$$K_1=n\omega+ t^d+ n^s \omega^{2s+2} + (\omega-t)^d.$$
Since $1\le t\le \omega-1$, it follows that $t^d+(\omega-t)^d\le (\omega-1)^d+1$ (since $d\ge 1$),
and so 
$$K_1=n\omega+ t^d+ n^s \omega^{2s+2} + (\omega-t)^d\le n\omega+n^s \omega^{2s+2}+ (\omega-1)^d+ 1= \omega^{s+3}+\omega^{s^2+4s+2}+(\omega-1)^d+1.$$
Hence $M^*$ is $(K,D)$-colourable where
$$K =\omega^{s+3}+\omega^{s^2+4s+2}+(\omega-1)^d+1.$$
This proves \ref{largechi}.~\bbox

From \ref{smallchi} and  \ref{largechi}, with $m=(s+1)\omega^s$, we deduce:
\begin{thm}\label{nbrchi}
Suppose that $d\ge 1$ and $D\ge 0$ have the property that every $H_s$-free graph $G'$ with $\omega(G')<\omega(G)$ is
$(\omega(G')^d,D)$-colourable; and $c\ge 2s$ satisfies \ref{Kst} when $H=H_s$. Assume that $\omega^3>\omega+s/7$
 and $\omega\ge 200$.
Then $V(\mathcal{L})\cup N(\mathcal{L})$ is $(K,D)$-colourable where 
$$K=(\omega-1)^d+\omega^{(c+1)(s+7)}.$$
\end{thm}
\Proof $V(\mathcal{L})$ has cardinality at most $14\omega^{s+7}$ from the definition of an $s$-template and since $k+|L_0|\le \omega$
by \ref{clique}. Every vertex in $N(\mathcal{L})$ is pendant or dense or pure with respect to $\mathcal{L}$, by \ref{types}. 
At most $14^{s+2}\omega^{s^2+9s+14}$ are pendant, by \ref{pendant}, and 
the chromatic number of the set of all dense vertices is at most $(14\omega^{s+6})^{c+1}\omega$ by \ref{semidense}.
By \ref{smallchi} with $m= (s+1)\omega^s$, and \ref{largechi}, the set of all pure vertices is
$(K_1,D)$-colourable where
$$K_1 =(\omega-1)^d+\omega^{s^2+4s+2}+\omega^{s+3}+1+ (2s+2)\omega^{s+1}+2\omega^{s+2}.$$
Adding, we deduce that $V(\mathcal{L})\cup N(\mathcal{L})$ is $(K,D)$-colourable
where 
$$K=14\omega^{s+7}+14^{s+2}\omega^{s^2+9s+14}+(14\omega^{s+6})^{c+1}\omega+(\omega-1)^d+\omega^{s^2+4s+2}+\omega^{s+3}+1+ (2s+2)\omega^{s+1}+2\omega^{s+2}.$$
Since $K\le (\omega-1)^d+ \omega^{(c+1)(s+7)}$, this proves \ref{nbrchi}.~\bbox

\section{Proof of the main theorem}

In this section we prove \ref{mainthm}.
If $\mathcal{L}=(L_0\ll L_k)$ is an optimal $s$-template in $G$, we define $Z(\mathcal{L})$ to be the union of $L_0$ and the set of all 
vertices $v\in N(\mathcal{L})$ such that for $1\le i\le k$,
$v$ has at most $\omega^{s+2}/4$
non-neighbours in $L_i$. 
Let $Y(\mathcal{L})= (V(\mathcal{L})\cup N(\mathcal{L}))\setminus Z(\mathcal{L})$.
First we need:
\begin{thm}\label{outnbrs}
Let $s\ge 1$ be an integer, and let $G$ be $H_s$-free, with $\omega(G)\ge 15$ and $\omega^2>s+1$. Let $\mathcal{L}$ be an optimal $s$-template in $G$.
Then every vertex in $Y(\mathcal{L})$ has at most $\omega^{s+7}$ neighbours in $V(G)\setminus Y(\mathcal{L})$.
\end{thm}
\Proof
Let $v\in Y(\mathcal{L})$, and suppose that $v$ has more than $\omega^{s+7}$ neighbours in $V(G)\setminus Y(\mathcal{L})$. 
By \ref{clique} and \ref{dense}, $|Z(\mathcal{L})|\le \omega+ 14\omega^{s+6}$, and so $v$ has at least $\omega^s$
neighbours in $V(G)\setminus (V(\mathcal{L})\cup N(\mathcal{L}))$, since
$$ \omega+ 14\omega^{s+6}+\omega^s<\omega^{s+7}$$
(because $\omega\ge 15$). By \ref{Ramsey}, there is a stable set $S\subseteq V(G)\setminus (V(\mathcal{L})\cup N(\mathcal{L}))$
of neighbours of $v$, with $|S|=s$. Since $v\notin Z(\mathcal{L})$, it follows that $v\in N(\mathcal{L})$
and there exists $i\in \{1\ll k\}$
such that $v$ has more than $\omega^{s+2}/4$ non-neighbours in $L_i$. Since $k\ge 2$  by \ref{nottoosmall}, and $v$ has a neighbour 
in at least one of $L_1\ll L_k$ because $v\in N(\mathcal{L})$, and each $|L_j|\ge \omega^{s+5}\ge \omega^{s+2}/4$, we may choose
distinct $i,j\in \{1\ll k\}$ such that $v$ has a neighbour in $L_j$ and $v$ has at least $\omega^{s+2}/4$ non-neighbours in $L_i$,
and choose such a pair $i,j$ such that $v$ has as many non-neighbours in $L_i$ as possible. Let $B$ be the set of non-neighbours 
of $v$ in $L_i$. Thus $|B|\ge \omega^{s+2}/4$.
\\
\\
(1) {\em $v$ has at most $\omega^s+\omega^{s+3}$ non-neighbours in $L_j$.}
\\
\\
Let 
$u\in L_j$ be adjacent to $v$. If $u$ has at least $\omega^s$ neighbours in $B$, there is a stable set $T$ of such neighbours
with $|T|=s$, by \ref{Ramsey}, and then the subgraph induced on $S\cup T\cup \{u,v\}$ is isomorphic to $H_s$, a contradiction.
Thus $u$ has fewer than $\omega^s$ neighbours in $B$; but it has at most $\omega^{s+3}$ non-neighbours in $L_i$, and hence at most
that many in $B$, and so $|B|\le \omega^s+\omega^{s+3}$. Since $|L_i|\ge \omega^{s+5}$, $v$ has a neighbour in $L_i$;
and so from the choice of the pair $i,j$, $v$ has at most $\omega^s+\omega^{s+3}$ non-neighbours in $L_j$. This proves (1).

\bigskip

Since $|B|\ge \omega^{s+2}/4\ge \omega^s$, it includes a stable set $T$ of cardinality $s$, by \ref{Ramsey}. Each vertex in $S$
has at most $\omega^{s+3}$ non-neighbours in $L_j$, and $v$ has at most $\omega^s+\omega^{s+3}$ non-neighbours in $L_j$ by (1),
and since 
$$s\omega^{s+3}+\omega^s+\omega^{s+3}<\omega^{s+5}\le |L_j|,$$
(because $\omega^2>s+1$), it follows that some vertex $u\in L_j$ is adjacent to every vertex in $T\cup \{v\}$. But then 
the subgraph induced on $S\cup T\cup \{u,v\}$ is isomorphic to $H_s$, a contradiction. This proves \ref{outnbrs}.~\bbox

Let $A\subseteq V(G)$, 
and let $\mathcal{L}=(L_0\ll L_k)$ be an $s$-template of $G$ with $V(\mathcal{L})\subseteq A$. We make the following definitions: 
\begin{itemize}
\item $Z_A(\mathcal{L})$ is the union of $L_0$ and the set of vertices in $A\setminus V(\mathcal{L})$ that have 
at most  $\omega^{s+2}/4$
non-neighbours in $L_i$ for each $i\in \{1\ll k\}$; 
\item $N_A(\mathcal{L})$ is the set of vertices in $A\setminus V(\mathcal{L})$ with a neighbour in 
$L_1\cupcup L_k$; and 
\item 
$Y_A(\mathcal{L})=(V(\mathcal{L})\cup N_A(\mathcal{L}))\setminus Z_A(\mathcal{L})$.
\end{itemize}

\begin{thm}\label{local}
Let $s\ge 1$ be an integer, and let $G$ be $H_s$-free. Let $\omega=\omega(G)$, and let $A\subseteq V(G)$, 
such that there is an $s$-template $\mathcal{L}$ of $G$ with $V(\mathcal{L})\subseteq A$ and with value at least $28\omega^{s+6}$. 
Let $\mathcal{L}=(L_0\ll L_k)$ be such an $s$-template with maximum value. 
\begin{itemize}
\item Suppose that $d\ge 1$ and $D\ge 0$ have the property that every $H_s$-free graph $G'$ with $\omega(G')<\omega$ is
$(\omega(G')^d,D)$-colourable; and $c\ge 2s$ satisfies \ref{Kst} when $H=H_s$; and $\omega^3>\omega+s/7$; 
 and $\omega\ge 200$.
Then $V(\mathcal{L})\cup N_A(\mathcal{L})$ is $((\omega-1)^d+\omega^{(c+1)(s+7)},D)$-colourable.
\item 
If  $\omega\ge 15$ and $\omega^2>s+1$, then every vertex in $Y_A(\mathcal{L})$ has at most $\omega^{s+7}$ neighbours in $A\setminus Y_A(\mathcal{L})$.
\end{itemize}
\end{thm}
\Proof
$\mathcal{L}$ is not necessarily an optimal $s$-template of $G$, since it is constrained to have vertex set included in $A$; and it is
not necessarily an optimal $s$-template of $G[A]$, since perhaps $\omega(G[A])<\omega(G)$ and then the conditions that define an
$s$-template of $G$
are different from those that define an $s$-template of $G[A]$. Nevertheless, we can apply \ref{nbrchi} and \ref{outnbrs} to $\mathcal{L}$
by the following trick. Let $G'$ be the disjoint union of $G[A]$ and a 
complete graph $K_{\omega(G)}$ with vertex set $B$ say. Then $\omega(G')=\omega(G)$, and
since every optimal $s$-template of $G'$ induces a connected subgraph with more than $\omega(G)$ vertices (because its value is at least
$28\omega(G)^6$), it contains no vertex of $B$ and so is an $s$-template of $G[A]$. This proves that $\mathcal{L}$
is an optimal $s$-template of $G'$. The claims of the theorem follow by applying \ref{nbrchi} and \ref{outnbrs} to $\mathcal{L}$ and $G'$. This proves \ref{local}.~\bbox

Now we prove \ref{mainthm}, which we restate in a strengthened form:
\begin{thm}\label{mainthm2}
For every integer $s\ge 1$, there exists $d\ge 0$ such that if $G$ is $H_s$-free, then $G$ is $(\omega^d,\omega^{s+8})$-colourable,
and hence has chromatic number at most $\omega^d(\omega^{s+8}+1)$.
\end{thm}
\Proof Choose $c\ge 2s$ satisfying \ref{Kst} with $H=H_s$. It follows from the main theorem of \cite{kierstead} that there is a function $f$ such that
$\chi(G)\le f(\omega(G))$ for every $H_s$-free graph $G$; and so by choosing $d$ sufficiently large we may arrange that
$\chi(G)\le \omega(G)^d$ for every $H_s$-graph $G$ with $\omega(G)$ at most $\max(200, (s+1)^{1/2})$. Let us also choose $d$ so 
large that $d\ge (c+1)(s+7)+1$. We claim that $d$ satisfies \ref{mainthm2}. The proof is by induction on $\omega=\omega(G)$.
If $\omega\le \max(200, (s+1)^{1/2})$ the claim is true, so we may assume that $\omega>\max(200, (s+1)^{1/2})$. Consequently
$\omega^3>\omega+s/7$ and $\omega^2>s+1$, and we can apply \ref{local}.

Let $V_0=V(G)$, and choose $n$ maximum such that there is a sequence $\mathcal{L}_i\;(1\le i\le n)$ of $s$-templates of $G$ and a sequence
$V_0\supseteq V_1\supseteq \cdots \supseteq V_{n}$ of subsets of $V(G)$, with the following property. For $1\le i\le n$,
$\mathcal{L}_i$ is an $s$-template of $G$ with value at least $28\omega^{s+6}$, with $V(\mathcal{L}_i)\subseteq V_{i-1}$, chosen with maximum value
among all such $s$-templates; and
$V_i=V_{i-1}\setminus Y_{V_{i-1}}(\mathcal{L}_i)$, in the notation of \ref{local}.
For $1\le i\le n$, let $Y_i=Y_{V_{i-1}}(\mathcal{L}_i)$, and let $Y_{n+1}=V_n$. We observe:
\begin{itemize}
\item The sets $Y_1\ll Y_{n+1}$ are pairwise disjoint and have union $V(G)$. 
\item For $1\le i\le n$, $G[Y_i]$ is $(\omega^d, \omega^{s+7}(\omega-1))$-colourable. To see this, observe that 
from the inductive hypothesis, every $H_s$-free graph $G'$ with $\omega(G')<\omega$ is $(\omega(G')^d,\omega(G')^{s+8})$-colourable,
and hence $(\omega(G')^d,D)$-colourable where $D=\omega^{s+7}(\omega-1)$, since $\omega(G')<\omega$. 
From the first statement of \ref{local}, for $1\le i\le n$, $Y_i$ is 
$((\omega-1)^d+\omega^{(c+1)(s+7)},\omega^{s+7}(\omega-1))$-colourable and hence $(\omega^d, \omega^{s+7}(\omega-1))$-colourable,
since $(\omega-1)^d+\omega^{(c+1)(s+7)}\le \omega^d$ (because $d\ge (c+1)(s+7)+1$).
\item $G[Y_{n+1}]$ is $(\omega^d, \omega^{s+7}(\omega-1))$-colourable. To see this, observe that 
the maximality of $n$ implies that there is no $s$-template $\mathcal{L}$ of $G$ of value at least $28\omega^{s+6}$ and with 
$V(\mathcal{L})\subseteq V_n$, and so $\chi(V_n)\le (14\omega^{s+6})^c$ by \ref{Kst}. Consequently $V_n=Y_{n+1}$ is
$(\omega^d, \omega^{s+7}(\omega-1))$-colourable, since $ (14\omega^{s+6})^c\le \omega^d$.
\item For $1\le i\le n$, every vertex in $Y_i$ has at most $\omega^{s+7}$
neighbours in $Y_{i+1}\cupcup Y_{n+1}$.
This follows from 
the second statement of \ref{local}.
\end{itemize}
By \ref{chaining}, $G$ is $(\omega^d, \omega^{s+7}(\omega-1)+\omega^{s+7})$-colourable.
This proves \ref{mainthm2}.

\end{document}